\documentclass{article}
\usepackage{graphicx}
\usepackage{float}
\usepackage{topcapt}

\usepackage{lipsum}

\usepackage{amsmath}
\usepackage{url}
\usepackage[T2A]{fontenc}
\usepackage[russian,english,british]{babel}
\usepackage{longtable}
\usepackage{fancyvrb}
\usepackage{tocbibind}

\begin{document}

\title{"The Gibbon of Math History". Who Invented the St. Petersburg Paradox? Khinchin's resolution.}
\author{Valerii Salov}
\date{}
\maketitle

\begin{abstract}
A sentence from Carl Boyer's \textit{A History of Mathematics} can be interpreted so that the full brothers Nicolaus II (02/06/1695 - 07/31/1726) and Daniel Bernoulli (02/08/1700 - 03/17/1782) are the authors of the St. Petersburg paradox. The paradox was formulated by their cousin Nicolaus I Bernoulli (10/21/1687 - 11/29/1759). The author did not find evidences that Nicolaus II and Daniel Bernoulli discussed the paradox. The key articles on the topic and its history from Karl Menger and Paul Samuelson, and recent papers presenting the time resolution miss Alexandr Khinchin's resolution of the paradox in his paper "On Petersburg game." \textit{Matematicheskii sbornik}, Volume 32, No 2, 1925, pp. 330 - 341. The C++ program khinchin.cpp simulates conditions of two Khinchin's theorems and confirms his results providing concrete empirical dependencies of the frequencies corresponding to geometric and arithmetic mean payments on the number of Petersburg games.
\end{abstract}

\section{Introduction}

Carl Benjamin Boyer (11/02/1906 - 04/26/1976) has made a great contribution to the history of mathematics. David Foster Wallace names him "the Gibbon of math history" \cite[p. 8]{wallace2010}. He writes \textit{"Boyer is joined at the top of the math-history food chain only by Prof. Morris Kline. Boyer's and Kline's major works are respectively A History of Mathematics and Mathematical Thought from Ancient to Modern Times"}. The first edition \cite{boyer1968} is the only one published during Boyer's life. The second \cite{boyer1991} and third \cite{merzbach2011} editions have been revised by Uta Merzbach. Comparison with Edward Gibbon (09/08/1737 - 01/16/1794), a historian who's \textit{"The History of the Decline and Fall of the Roman Empire"} was published in six volumes between February 1776 and May 1788, emphasizes the fundamental character of Boyer's work. The time resolution of the St. Petersburg paradox \cite{peters2011}, \cite{peters2011b} and criticism of the proposal \cite{varma2013} have inspired the author to review the history of the topic. He has found the following sentence in \cite[p. 463]{boyer1968}, \cite[p. 423]{boyer1991}, \cite[p. 397]{merzbach2011}

\textit{"When Daniel Bernoulli went to St. Petersburg in 1725, his older brother also was called there as a professor of mathematics; in the discussions of the two men there arose a problem that has come to be known as the Petersburg paradox, probably because it first appeared in the Commentarii of the Academy there.$^6$"}

The top index 6 and the footnote are present only in \cite[p. 463]{boyer1968}

\textit{"The full title is Commentarii Academiae Scientiarum Imperialis Petropolitanae. This journal contains many articles by the younger Bernoullis and their Swiss colleague, Euler."}

This sentence creates an impression that inventors of the St. Petersburg paradox are the older and younger full brothers Nicolaus and Daniel Bernoulli. Another interpretation might be: both knew about the paradox, could discuss it rapidly since worked at one academy, and developed the topic to such an extent that it has become famous. There is some sense in the sentence that Daniel went first and then Nicolaus was called to the St. Petersburg Academy. Let us investigate.

\section{Nicolaus I, Nicolaus II, and Daniel Bernoulli}

Boyer writes \cite[p. 455]{boyer1968} \textit{"No family in the history of mathematics has produced as many celebrated mathematicians as did the Bernoulli family ..."}. He creates \textit{"The mathematical Bernoullis: a genealogical chart"} \cite[p. 456]{boyer1968}, \cite[p. 416]{boyer1991}, where the mathematicians of our primary interest are Nicolaus II (1687 - 1759), Nicolaus III (1697 - 1726), and Daniel I (1700 - 1782). The chart contains 13 names: four Nicolauses, three Jeans, two Daniels, two Jacqueses, one Christoph, and one Jean Gustav. A way to distinguish the members is needed. Another complication is different spelling of the first names occurred because several members lived and worked in different countries. The names cited in Swiss, Spanish (S), Anglicized (A), German (G), and Russian (R) forms are: Jacques (S) - James (A) - Jacob (G) - Yakov (R), Jean (S) - John (A) - Johann (G) - Iogann, Ivan (R), Daniel (S, A) - Daniil (R), Nicolaus (S) - Nicolas (A) - Nicolai (R). Using the Roman numerals is shorter than the years of life. In this paper, the author applies the numbers accepted in \cite{yushkevich1972}, \cite{fleckenstein2008}, \cite{fleckenstein2008b}, \cite{straub2008}, and the two Internet sources \url{http://www-history.mcs.st-and.ac.uk/BiogIndex.html} and \textit{Wikipedia}: Nicolaus I (1687 - 1759), Nicolaus II (1695 - 1726), Daniel (1700 - 1782). The former Internet source is used for citing the dates of birthday and death. The dates were checked when possible with Complete Dictionary of Scientific Biography \cite{straub2008}. All dates are of the new style unless the old style is labeled explicitly. These numbers are used by Jay Goldman \cite[p. 25]{goldman1998}: \textit{"Nicolas II and Daniel Bernoulli (Johann's sons) went to Russia's newly formed St. Petersburg Academy in 1725"} and in the \textit{scientific biography} of Daniel Bernoulli prepared by Ashot Tigranovich Grigoryan and Boris Demyanovich Kovalev \cite{grigoryan1981}. On page 14 of the latter we find a genealogical tree with 15 members of Bernoulli family. Nicolaus I is a cousin of Nicolaus II and Daniel.

The three volumes \cite{yushkevich1972} is a result of a collective work headed by Adolf Pavlovich Yushkevich (07/15/1906 - 07/17/1993). The third volume is devoted to the XVIII century. The most relevant for our consideration Chapter Four "Probability Theory" pp. 126 - 152 of the third volume is written by Oscar Borisovich Sheinin (1925 - ) and Leonid Efimovich Maistrov (1920 - 1982). It contains Sections "From Y. Bernoulli to Moivre" pp. 126 - 128 and "Works of D. Bernoulli", pp. 140 - 144. Only in the third volume the names of Nicolaus I, Nicolaus II, and Daniel Bernoulli are referenced 20, 13, and 64 times.

On page 18 \cite{yushkevich1972} there is a photocopy of the title-page of Tomus I of "Commentarii Acedemiae Scientiarum Imperialis Petropolitanae" issued in 1726, where the relevant to the discussion memoir of Daniel Bernoulli was published in 1738. Page 20 presents a photocopy of the beginning list of an advertisement about the lectures to be read in 1726 since January 24, 1726 (old style) on Mondays, Wednesdays, Thursdays, and Saturdays. The first three names in the list are Daniil Bernoulli (hour from 7 to 8), Teofil Zigfrid Baer, Nicolai Bernoulli (hour from 8 to 9). With a reference to \textit{L.G. du Pasquier. Leonard Euler et ses amis. Paris, 1927, p. 9} there is a citation of the words of their father Johann Bernoulli (06/06/1667 - 01/01/1748) (see biography in \cite{fellmann2008}) on pages 21 - 22 (author's translation from Russian): \textit{"It is better to be patient for a while with severe climate of the country of ices, in which muses are greeted, than to die from hunger in a country with mild climate in which the muses are offended and despised."} Petr Petrovich Pekarskii (05/31/1827 - 07/24/1872) writes the biography of Nicolaus II Bernoulli \cite[pp. 95 - 98]{pekarskii1870}. The abstract lists primary sources including biography composed by Christian Goldbach (03/18/1690 - 11/20/1764) soon after death of Nicolaus II and published in Commentarii II 482, II 266 - 270. The abstract mentions pages 229, 230, 239 - 291, 294, 295, 298, 299 describing friendship between two brothers Nicolaus and Daniel Bernoulli. The following author's translation from Russian repeats the old style dates:

\textit{"Both brothers arrived in Petersburg on 27 October 1725. ... After an eight months stay in Petersburg, Nicolaus Bernoulli was taken ill and died on 29 July 1726 from abscess of viscera as it was diagnosed after corpse lancing. Three days later there was a ceremonial meeting of the academy visited by empress Ekaterina. Having learned about the death of the academician, she called to herself his brother and consoled him in mercy expressions. Manuscripts of Nicolaus Bernoulli were given to brother Daniel and the following two papers were published in Commentarii Academiae petropolitanae: I. 121 - 126 De motu corporum ex percussione; 198 - 207 Analysis aequationum quarundam differentialium"}.

The inviting letters, salaries, contracts signed between the brothers and Academy, temporary confusion with the first names are described in \cite[pp. 43 - 46]{grigoryan1981}. Nicolas II and Daniel got 1000 and 800 rubles per year. The starting academic salary 1000 rubles was the highest in 1725. The author has found a typo in the cited date January 28, 1728 \cite[p. 100]{pekarskii1870} of the letter sent by Daniel to Goldbach. The year must be 1725. Daniel describes him as a young 25 years old man and asks Godlbach, who was in January 1725 in St Petersburg, to let know about him to Blumentrost and Golovkin. He just got a letter from Nicolaus II: \textit{"... he, from the true brothers friendship, decides not let me go alone to Moskoviya and agrees to sacrifice his benefits (a chair brining 150 Luidors) and accompany me. I believe, it would be easy for both of us to find positions in St. Petersburg especially as there is nothing more spacious  than mathematics."} (author's translation from Russian; Moskoviya means Russia).

Nicolaus I, Johann's nephew, had not been in St. Petersburg. While Boyer's sentence does not name the "older brother", Boyer talks about discussions between Daniel and his \textit{full brother} Nicolaus II. Both arrived in St. Petersburg together and regrettably had no more than nine months for discussions.

\section{Specimen Theoriae Novae de Mensura Sortis. Commentarii (1730 - 1731) 1738}

From seven memoirs of Daniel Bernoulli on probability theory six were published in Commentarii \cite[p.140]{yushkevich1972}. Louise Sommer has translated to English the one in the title \cite{bernoulli1954}. Bernoulli was writing in Latin and citing the letter of Gabriel Cramer (07/31/1704 - 01/04/1752) in French. Translation contains four footnotes made by Karl Menger. Bernoulli confirms who is the author of the problem known today as the St. Petersburg paradox or game \cite[31]{bernoulli1954}

\textit{"My most honorable cousin the celebrated Nicolaus Bernoulli, Professor utriusque iuris at the University of Basle, once submitted five problems to the highly distinguished mathematician Montmort. These problems are reproduced in the work L'analyse sur les jeux de hazard de M. de Montmort, p. 402. The last of these problems runs as follows: Peter tosses a coin and continues to do so until it should land "heads" when it comes to the ground. He agrees to give Paul one ducat if he gets "heads" on the very first throw, two ducats if he gets it on the second, four if on the third, eight of on the fourth, and so on, so that with each additional throw the number of ducats he must to pay is doubled. Suppose we seek to determine the value of Paul's expectation. My aforementioned cousin discussed this problem in a letter to me asking for my opinion."}

The letter from Nicolaus I to Pierre R\'{e}mond de Montmort (10/27/1678 - 10/07/1719) mentioned by Gabriel Cramer writing to Nicolaus I  \cite[p. 33]{bernoulli1954}: \textit{... I believe that I have solved the extraordinary problem which you submitted to M. de Montmort, in your letter on September 9, 1713, (problem 5, page 402)"} was appended to the second edition of \textit{Essai d'Analyse sur les Jeux de Hazard, Paris, 1713}, see \cite[p. 31, Translator's footnote 8]{bernoulli1954}, \cite{montmort1713}. The author have noticed that the publication year in e-rara reference is 1718. However, the title-list \cite[p. 8 of the photocopy]{montmort1713} has the Roman MDCCXIII year. This is 1713. Cramer's letter to Nicolaus I containing Cramer's resolution was sent in 1728 \cite[p. 33]{bernoulli1954}. Georges Louis Leclerc Comte de Buffon (09/07/1707 - 04/16/1788) recollects \cite[p. 75, Section XV]{buffon1777} that during his visit to Geneva in 1730 Cramer presented him this task, said that first time Nicolaus I Bernoulli suggested it to Montmort, and that the task could be found on pages 402 - 407 in Montmort's book. The author has noticed that the entire letter of Nicolaus I Bernoulli in Montmort's book occupies two pages 401 and 402 in the photocopy of the original \cite{montmort1713}. Buffon also mentions Daniel Bernoulli's resolution and some similarity of Cramer's and Bernoulli's proposals. Let us make the dated list of items:

\begin{enumerate}
\item[1)]
September 9, 1713. Letter of Nicolaus I Bernoulli to de Montmort describing the problem named later the St. Petersburg paradox or game \cite[p. 402, \textit{Cinqui\'{e}me Probl\'{e}me}]{montmort1713};
\item[2)]
1713. Second edition of "Essai d'Analyse sur les Jeux de Hazard", Paris publishing the letter from point 1 \cite[p. 401 - 402]{montmort1713};
\item[3)]
???? There was a letter of Nicolaus I to Daniel. Daniel Bernoulli: \textit{"My aforementioned cousin discussed this problem in a letter to me asking for my opinion"} \cite[p. 31]{bernoulli1954};
\item[4)]
1728. Solution of the problem sent to Nicolaus I Bernoulli by Cramer in a letter \cite[p. 33]{bernoulli1954}. Daniel Bernoulli knew about it in 1732;
\item[5)]
1730. A verbal discussion between Cramer and Buffon in Geneva about Nicolaus's I problem and Cramer's resolution \cite[p. 75]{buffon1777};
\item[6)]
1730 - 1731. Daniel Bernoulli solves the problem without knowing about the Cramer's solution, differently, and \textit{reads it to the Society} \cite[p. 33, Translator's footnote 11 clarifies that \textit{Society} is the Imperial Academy of Sciences in Petersburg]{bernoulli1954}.
\item[7)]
1731 - 1732. There was a letter of Daniel to Nicolaus I: \textit{"After having read this paper to the Society I sent a copy to the aforementioned Mr. Nicolaus Bernoulli, to obtain his opinion of my proposed solution to the difficulty he had indicated."} \cite[p. 33]{bernoulli1954}.
\item[8)]
1732. A letter of Nicolaus I to Daniel. Daniel Bernoulli: \textit{"In a letter to me written in 1732 he declared that he was in no way dissatisfied with my proposition ... Then this distinguished scholar informed me that the celebrated mathematician, Cramer, had developed a theory on the same subject several years before I produced my paper. ... Cramer himself first described his theory in his letter of 1728 to my cousin."} \cite[p. 33]{bernoulli1954};
\item[9)]
1738. Publishing "Specimen Theoriae Novae de Mensura Sortis" \cite{bernoulli1954}.
\end{enumerate}

If the letter in point three was sent to Daniel after July 1726, then brothers did not discuss the problem. The letter could clarify when Daniel Bernoulli was involved. The second edition of Montmort's book \cite{montmort1713} existed for 13 years but Daniel knew about the task from his cousin's letter. Could Nicolaus II know the task from \cite{montmort1713}? Even if the latter is true, then Daniel's reference to cousin's letter indicates that Nicolaus II unlikely said about it to Daniel. Was the mentioned letter sent prior July 1726?

The eight works of Nicolaus II Bernoulli \cite{fleckenstein2008b} including the two mentioned by Pekarskii \cite{pekarskii1870} given to Daniel Bernoulli after the death of his brother are not on the probability theory.

Daniel Bernoulli's thesis \textit{On Breath} was defended in 1720 \cite[p. 427]{bernoulli1959}. The \textit{Dissertatio inauguralis physico-medica de respiratione (Basel, 1921)} \cite[p. 44, reference 1]{straub2008} is far from the probability theory. By own admission, mathematics attracted him due to the \textit{"family members example"} and \textit{"own soul propensity"} \cite[p. 428]{bernoulli1959}. He says about the period 1720 - 1723. Year 1724 marks his first mathematical publication \textit{Mathematical Exercises} \cite[p. 428]{bernoulli1959}, \cite[p. 44, reference 1]{straub2008}: \textit{Exercitationes quaedam mathematicae (Venice, 1724)}, \textit{"which attracted so much attention that he was called to the St. Petersburg Academy"} \cite[p. 36]{straub2008}. Hans Straub summarizes the four topics of this treatise \cite[p. 37]{straub2008}: \textit{"... the game of faro, the outflow of water from the openings of containers, Riccati's differential equation, and the lunulae (figures bounded by two circular arcs)"}. A detailed review of the book, scientific, and cultural background of the epoch we find in \cite[pp. 22 - 41, Chapter 2 Debuts]{grigoryan1981}. It marks erudition of young Daniel Bernoulli but does not show that his book discusses \textit{Cinqui\'{e}me Probl\'{e}me} from \cite[p. 402]{montmort1713}. Examining the photocopy \cite{bernoulli1724}, the author found several pages mentioning Montmort but unrelated to the Nicolaus's I paradox. 

Geography and time of travels of Nicolaus II \cite[p. 97]{pekarskii1870} and Daniel \cite[pp. 99 - 100]{pekarskii1870} Bernoulli during 1720 - 1725 described by Pekarskii do not create long intersections.

Nicolaus II and Daniel Bernoulli were not only brothers but friends \cite[p. 95]{pekarskii1870}. Daniel Bernoulli carefully cites contributions of his cousin and Cramer \cite{bernoulli1954}. His own contribution was indisputable and known since \textit{reading to the Society} in 1731. There were no reasons for him to be silent with respect to a Nicolaus's II contribution to the Nicolaus's I game, if such would exist in a useful noticeable form, even, as a verbal discussion.

Based on the reviewed sources, the author concludes:
\begin{enumerate}
\item[a)]
the inventor of the St. Petersburg paradox or game is Nicolaus I Bernoulli, who reported about it in the published letter dated by September 9, 1713; this was known prior Carl Boyer's publication \cite{boyer1968}, for instance, from \cite{bernoulli1954}; 
\item[b)]
no evidences confirming intensive discussions of the paradox between Nicolaus II and Daniel Bernoulli are found; both arrived in St. Petersburg on October 27, 1725 (old style);
\item[c)]
the date of the letter mentioned in item 3 could put more light, especially, if it was sent after July 1726; the letter could be sent by Nicolaus I to Daniel Bernoulli after July 1726 and prior Nicolaus I Bernoulli received Cramer's letter dated by 1728;
\item[d)]
in the date 28 January 1728 on page 100 in \cite{pekarskii1870}, the year must be 1725; this is the date of the letter sent by Daniel Bernoulli to Christian Goldbach;
\item[e)]
there is discrepancy in the publication year on the title page of \cite[p. 8 of photocopy]{montmort1713}, where the Roman MDCCXIII means 1713, and the e-rara.ch project \cite{montmort1713} dating this publication by 1718; the five years difference is important for better understanding of \textit{chances} that Nicolaus II Bernoulli knew about \textit{Cinqui\'{e}me Probl\'{e}me} from the book prior his arrival in St. Petersburg; it should be clarified or corrected;
\item[f)]
there is discrepancy in Buffon's recollection \cite[p. 75]{buffon1777} that Nicolaus's I problem can be found on pages 402 - 407 in Montmort's book and the photocopy of the available edition \cite{montmort1713} presenting Nicolaus's I letter on two pages 401 and 402; the problem is on page 402;
\item[g)]
the cited sentence found in the three editions of famous Carl Boyer's and Uta Merzbach's book should be changed in order to avoid confusion.
\end{enumerate}

\section{Khinchin's Resolution}

Vladimir Ivanovich Smirnov (06/10/1887 - 02/11/1974) writes a condensed and thorough article on Daniel Bernoulli \cite{smirnov1959}, where pages 461 - 465 are devoted to the St. Petersburg paradox, Daniel Bernoulli's resolution, and later works related to the problem. His consideration is finished by the words (author's translation from Russian): \textit{"Explanation of the "paradox" of the "Petersburg game" was given by A. Ya. Khinchin in his work "On Petersburg game" (Matematicheskii sbornik, 1925). In this work a question on some average winning for a big number of games is investigated."} \cite[p. 465]{smirnov1959}. The paper of Alexandr Yakovlevich Khinchin (07/19/1894 - 11/18/1959) is published in Russian together with abstract in German \cite{khinchin1925}.

At the beginning, he presents \textit{"a known Petersburg game"} in expressions close to Daniel Bernoulli's description and confirms that \textit{"for explanation of the "paradox" many considerations were proposed but we shall not stop on them here"} (author's translation from Russian is here and below). \textit{"Let us notice only that in this case, of course, no speech may go about any mathematical paradox but at most about that the mathematical expectation is not always adequate to those worldly-psychological representations, which it is commonly connected to. In the case of the Petersburg game, it is often pointed to that Petr in his expectation of winning, naturally, orients not on the mathematical expectation of winning in a particular game, which is difficult to account psychologically, but on some average winning during big number of games. Such understanding of psychological prerequisites of the "paradox" puts in front of us a certain mathematical task, which can be formulated as follows: Find such an estimate of the mean winning of Petr during a big number of games, that its probability would go to unit with infinite increasing the number of games. However, it makes sense to say, that the task will get a quite determined sense only after a certain notion of the mean winning will be exactly defined. In the current note, we shall consider in details the set problem in two of the most simple (and also the most important) cases, namely in assumption that the mean winning is defined as the geometric and arithmetic mean of particular games"} \cite[p. 330]{khinchin1925}. Accordingly, the paper formulates and proves two theorems.

The first theorem preceded by two Lemmas is formulated for the Peter's winning $a_n$ in $n$th game after introducing the denomination $\varrho_n = (a_1 a_2 \dots a_n)^{\frac{1}{n}}$.

\textit{"Theorem I. For any small $\delta > 0$ and $\eta > 0$, there exists the number $N = N(\delta, \eta)$, possessing the following property: with the probability exceeding $1 - \eta$, it is possible to expect, that for all $n > N$ we shall have $2 - \delta < \varrho_n < 2 + \delta$."} \cite[p. 332]{khinchin1925}.

Khinchin comments \cite[p. 333]{khinchin1925}: \textit{"In other words, with a probability arbitrary close to unit, it is possible to expect, that the geometric mean of Petr's winnings for sufficiently big number of games will be arbitrary close to 2 ducats"}.

The second theorem preceded by five Lemmas is formulated using the denomination $\sigma_n = \frac{1}{n}(a_1 + a_2 + \dots + a_n)$.

\textit{"Theorem II. For any small positive numbers $\delta > 0$ and $\eta > 0$, there exists the number $N = N(\delta, \eta)$, possessing the following property: with the probability exceeding $1 - \eta$ it is possible to expect that for all $n > N$ we shall have $1 - \delta < \frac{\log \sigma_n}{\log\log n} < 1 + \delta$"} \cite[p. 338]{khinchin1925}.

Khinchin comments \cite[p. 338]{khinchin1925}: \textit{"Otherwise speaking, with a probability arbitrary close to unit, it is possible to expect, that the order of growth of the arithmetic mean of Petr's winnings for the first $n$ games will be between $(\log n)^{1-\delta}$ and $(\log n)^{1+\delta}$"}.

\section{Uncited Results}

Karl Menger (01/13/1902 - 10/05/1985) has prepared a historical and mathematical review of the St. Petersburg game \cite{menger1934}. With his help we have the English translation of Daniel Bernoulli's work \cite{bernoulli1954} made by Louise Sommer. For many years Bernoulli's \textit{moral expectation}, Cramer's \textit{utility}, later progress and criticism of the \textit{utility theory}, and both articles remained in sight of leading economists Maurice Allais (05/31/1911 - 10/09/2010) \cite{allais1953}, Paul Samuelson (05/15/1915 - 12/13/2009) \cite{samuelson1960}, \cite{samuelson1977}, G\'{e}rard Debreu (07/04/1921 - 12/31/2004) \cite{debreu1951}, Kenneth Arrow (08/23/1921 - ) \cite{arrow1971}, Harry Markowitz (08/24/1927 - ) \cite{markowitz1952}, Daniel Kahneman (03/05/1934 - ) \cite{kahneman1979}, \cite{tversky1992} \cite{kahneman2002}, William Sharpe (06/16/1934 - ) \cite[Chapter 3]{sharpe2007}, Amos Tversky (03/16/1937 - 06/02/1996) \cite{kahneman1979}, \cite{tversky1992}. Borrowing the words of Yakov Sinai \cite[p. 5]{sinai2010}, \textit{"there are a number of glaring omissions in this text"}. However, Menger's paper does not reference Khinchin's work.

Second Samuelson's paper on the St. Petersburg paradox \cite{samuelson1977} is another milestone publication on its history and essence. It is silent with respect to Khintchin's article. It has a reference to the English translation of Maistrov's book. On pages 315 and 316 of the original \cite{maistrov1967}, we find the references 61 and 62 coinciding with \cite{smirnov1959}. Therefore, Samuelson was one step apart from Smirnov's reference to Khintchin.

Interest to the Petersburg game does not weaken. However, recent publications \cite{hayden2009}, \cite{peters2011}, \cite{peters2011b}, \cite{varma2013} do not mention \cite{khinchin1925}.

Khinchin is indifferent to that how the number of games $n$ is realized \cite[p. 339]{khinchin1925}: \textit{"We will assume again this state of affairs the carried out"}. If it is treated as the product $n = \lfloor \nu \times t \rfloor$, where $\nu$ is the frequency of games and $t$ is the play time, then the latter can be compared with a life duration by a person making decisions. The Petersburg game has a random duration, if the time of a coin trial is fixed. Under high frequency of trials and games (on a computer), if $n$ is realized, then its nature is irrelevant for Khinchin's theorems: choosing an \textit{ensemble} of \textit{independent} games played simultaneously or their time chain makes no difference, if $n$ is identical in both cases. The author believes that Khinchin's results obtained in 1925 are important for \cite{peters2011}. He "shyly" assumes that Theorem I on geometric mean, being better known, would be considered by John Kelly Jr. (05/24/1927 - 03/02/1985), when writing \cite{kelly1956}.

\section{Khinchin's Theorems and C++ Experiment}

Buffon made an experiment and found that in 2,048 games the payment to the winning side was 10,057 crowns \cite[p. 424]{boyer1991}. This is $\frac{10,057}{2,048} \approx 5$ crowns per game. Boyer cites the numbers exactly. Originals are \textit{"deux mille quarante-huit"} and \textit{"dix mille cinquante-fept"} \cite[p. 84, Section XVIII, paragraph 2]{buffon1777}. Buffon's units are \textit{"\'{e}cus"}. Samuelson discusses Buffon's views on the game \cite[pp. 39 - 41, \textit{Buffon's practical objections}]{samuelson1977} and a perspective of computer experiments.

One curious \textit{"Why do we believe theorems?"} can be redirected to Andrzej Pelc's \cite{pelc2011}. We do not need a computer program to convince us in \textit{Pythagorean theorem}. However, accurately measuring \textit{right triangles} drawn on an Earth surface, we can determine that $c^2 = a^2 + b^2 + O(\frac{1}{R^2})$ and a program can be useful. Buffon confirms absurdity of paying too many \'ecus for a game. We jump over 301 years (September 9, 1713 and 2013 are the paradox's birthday and 300th anniversary) to Bjarne Stroustrup's \textit{C++} \cite{stroustrup2000} and a computer simulation because Khinchin does not tell how $N$ depends on $\delta$ and $\eta$.

\textit{The C++ Standard Committee, JTC1/SC22/WG21} has \textit{"made a present to C++ programmers"} \cite[p. 78]{salov2013} adding Section \textit{26.5 Random number generation} to the \textit{ISO/IEC 14882:2011 Programming Language C++ draft}, where \textit{mt19937} and \textit{ mt19937\_64} are aliases of 32 and 64-bit uniform generator template class \textit{std::mersenne\_twister\_engine} representing Makoto Matsumoto's and Takuji Nishimura's invention \cite{matsumoto1998}. This will support our \textit{pseudo coin trials}.

Introducing the number $\eta > 0$, Khinchin applies it only in order to say that $N$ depends on both $\delta$ and $\eta$. The abstract in German does not mention $\eta$. No explicit estimate of $N$ based on given $\delta$ and $\eta$ is presented. Yuri Vasilevich Prokhorov's (12/15/1929 - 07/16/2013) estimate for the \textit{law of big numbers} \cite{prokhorov1983} is an illustration that such dependencies can be comprehensive. The probability inequality $P\{|\frac{\mu_n}{n} - p| \le \epsilon \} > 1 - \eta$, where $\mu_n$ is the number of successes in $n$ Bernoulli trials each with probability $p$, is satisfied for any two numbers $\epsilon > 0$ and $\eta > 0$, if $n_0 > \frac{1 + \epsilon}{\epsilon^2} \log \frac{1}{\eta} + \frac{1}{\epsilon}$. Let $\epsilon = 0.001$ and $\eta = 0.001$, then $n_0 > 6,915,663$.

\begin{figure}[h!]
  \centering
  \includegraphics[width=130mm]{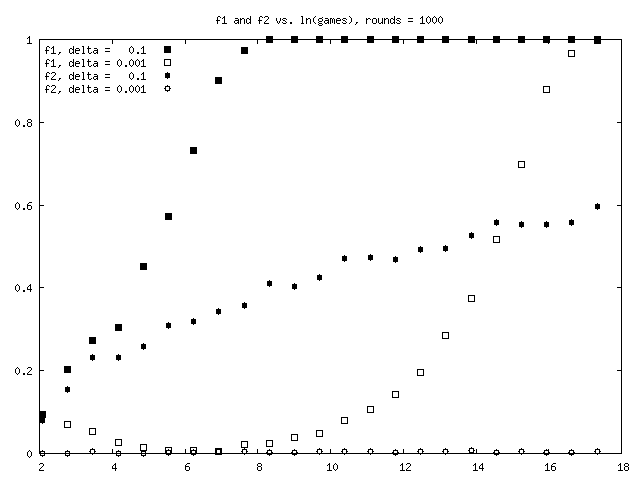}
  \caption[f1_f2_vs_g]
   {Dependence of experimentally observed frequencies $f1$ and $f2$ on natural logarithm of the number of games for different $\delta$. The frequencies correspond to Khinchin's geometric $\varrho_{\textrm{games}}$ and arithmetic $\sigma_{\textrm{games}}$ mean payments observed in the St. Petersburg game.}
  \label{f1_f2_vs_g}
\end{figure}

In C++ khinchin.cpp, Appendix, a \textit{croupier} tosses a fair, $p = \frac{1}{2}$, \textit{coin}. This is repeated until \textit{HEAD} is drawn. The number of \textit{tails} is counted. This completes one game with the payment $2^\textrm{tails}$. Croupier makes \textit{games}-number of repetitions and the geometric and arithmetic mean payments per game are computed. In theory, the \textit{while-loop} can iterate endlessly. The parameter $\delta$ is needed in order to check that $|\varrho_{\textrm{games}} - 2| < \delta$ and $[\log(\textrm{games})]^{1 - \delta} < \sigma_{\textrm{games}} < [\log(\textrm{games})]^{1 + \delta}$. Several \textit{rounds} are done in order to evaluate how many times both inequalities involving $\varrho_{\textrm{games}}$ and $\sigma_{\textrm{games}}$ take place. Dividing both counts by rounds, the program reports two frequencies $f1$ and $f2$ corresponding to two probabilities in Khinchin's theorems. The C++ built-in data types \textit{double}, \textit{unsigned int}, \textit{size\_t} are selected to increase the length of runs and accuracy of computations. The parameter $\eta$ might be useful, if the program should find by increasing the number of games and rounds the smallest $N$s for which $f1$ and $f2$ exceed $1 - \eta$. The number of rounds affects the accuracy of $f1$ and $f2$ estimated as ratios of integers. The task is omitted because convergence of $\sigma_{\textrm{games}}$ is found slow for $\textrm{games} \in [8, \; 33,554,432]$ and $\delta < 0.1$. In contrast, convergence of $\varrho_{\textrm{games}}$ is faster and $f1$ approaching 1 is observed routinely with this program for $\delta = 0.01$. Figure \ref{f1_f2_vs_g} illustrates typical dependencies.

Let us notice that the sample mean and standard deviation in 10 rounds presented in Appendix for A are 7.2 and 4.5, and $\log(2048) \approx 7.6$. Buffon's mean was 4.9 in a single round of 2048 games.

The Khichin's paper relates to his work on the \textit{law of iterated logarithm} partly summarized in \cite{khinchin1932}. The latter results have been generalized by Andrey Nikolaevich Kolmogorov (04/25/1903 - 10/20/1987) \cite{kolmogorov1929}.

\section{Appendix}

The C++ program khinchin.cpp is a \textit{command line application}. After compilation and linking, being run without arguments, it outputs help
\begin{verbatim}
C:\bin>khinchin
Usage: khinchin games delta [rounds = 100] [seed = 1234567] [details = no]
\end{verbatim}
The arguments games and delta corresponding to Khinchin's $n$ and $\delta$ are mandatory. The arguments rounds, seed, and details are optional. A run for Buffon's number of games 2048 using default rounds, seed, and details is
\begin{verbatim}
C:\bin>khinchin 2048 0.05
g = 2048 d = 0.05 r = 100 f1 = 0.76 f2 = 0.16 s = 1234567
\end{verbatim}
Additional output containing G geometric and A arithmetic means is generated after adding any text as the sixth argument counting from the program name
\begin{verbatim}
C:\bin>khinchin 2048 0.05 10 1234567 yes
G = 1.94263 A = 6.97852
G = 2.02246 A = 5.64502
G = 2.09634 A = 14.1782
G = 1.95847 A = 4.32813
G = 1.98114 A = 6.52002
G = 2.07376 A = 4.80371
G = 1.97645 A = 5.91211
G = 2.01563 A = 10.0278
G = 2.01631 A = 15.6226
G = 1.99055 A = 5.01318
g = 2048 d = 0.05 r = 10 f1 = 0.7 f2 = 0.1 s = 1234567
\end{verbatim}
The program rejects meaningless input
\begin{verbatim}
C:\bin>khinchin oh ah oi
games = oh, delta = ah, rounds = oi must be > 0 and games > 2
\end{verbatim}
The \textit{standard output} contains two fixed numbers of \textit{fields}. This is friendly for \textit{batch} and \textit{text processing} using AWK \cite{aho1988}, sed \cite{dougherty1997}, and \textit{graphics processing} using gnuplot \cite{gnuplot} and computing statistics of means.

\begin{verbatim}
#include <iostream>     /* Standard C++ cout, cerr */
#include <string>       /* Standard C++ string */
#include <cstdlib>      /* Standard C++ atoi */
#include <stdexcept>    /* Standard C++ invalid_argument */
#include <random>       /* Standard C++ random */
#include <cmath>        /* Standard C++ log, pow */
using namespace std;    /* Standard C++ namespace */

int main(int argc, char* argv[])
{
    try {
        const int       ROUNDS = 100;   // default number of rounds
        const int       SEED = 1234567; // default seed
        if(argc < 3) {
            cout    << "Usage: " << argv[0] << " games delta [rounds = "
                    << ROUNDS << "] [seed = " << SEED
                    << "] [details = no]" << endl;
            return	0;
        }
        const size_t    games = atoi(argv[1]);
        const double    delta = atof(argv[2]);
        const size_t    rounds = argc > 3 ? atoi(argv[3]) : ROUNDS;
        const int       seed = argc > 4 ? atoi(argv[4]) : SEED;
        if(games < 3 || delta <= 0.0 || rounds < 1)
            throw   invalid_argument(string("games = ") + argv[1] +
                ", delta = " + argv[2] + ", rounds = " + argv[3] +
                " must be > 0 and games > 2");
        const double    lnGames = log(static_cast<double>(games));
        const double    aMeanLow = pow(lnGames, 1.0 - delta);
        const double    aMeanHigh = pow(lnGames, 1.0 + delta);
        double          f1 = 0.0;   // frequency for Theorem I
        double          f2 = 0.0;   // frequency for Theorem II
        enum {TAIL, HEAD};
        uniform_int_distribution<int>   coin(TAIL, HEAD);
        mt19937                         croupier(seed);
        for(size_t r = 0; r < rounds; r++) {
            double  gMean = 0.0;    // geometric mean
            double  aMean = 0.0;    // arithmetic mean
            for(size_t g = 0; g < games; g++) {
                double  tails = 0.0;
                while(coin(croupier) == TAIL) tails++;
                gMean += tails; aMean += pow(2.0, tails);
            }
            gMean = exp(gMean * log(2.0) / games); aMean /= games;
            if(argc > 5)
                cout    << "G = " << gMean << " A = " << aMean << endl;
            if(fabs(gMean - 2.0) < delta) f1++;
            if(aMeanLow < aMean && aMean < aMeanHigh) f2++;
        }
        f1 /= rounds; f2 /= rounds;
        cout    << "g = " << static_cast<unsigned int>(games) << " d = "
                << delta << " r = " << static_cast<unsigned int>(rounds)
                << " f1 = " << f1 << " f2 = " << f2 << " s = " << seed
                << endl;
    }
    catch(const exception& e) {
        cerr    << e.what() << endl;
        return  -1;
    }
    catch(...) {
        cerr    << "Unknown exception" << endl;
        return  -2;
    }
    return  0;
}
\end{verbatim}
Compilation and linking on Windows using the Microsoft Visual C++ beginning from 2012 for a better support of C++ 11 is achieved with the command
\begin{verbatim}
C:\projects\consapp>cl /O2 /EHsc /Fekhinchin.exe khinchin.cpp
\end{verbatim}
The same in a UNIX environment using a GCC supporting C++ 11 is done as
\begin{verbatim}
$ c++ -O2 -std=c++11 -o khinchin.exe khinchin.cpp
\end{verbatim}

{

\bigskip

\noindent\textbf{Valerii Salov} received his M.S. from the Moscow State University, Department of Chemistry in 1982 and his Ph.D. from the Academy of Sciences of the USSR, Vernadski Institute of Geochemistry and Analytical Chemistry in 1987.  He is the author of the articles on analytical, computational, and physical chemistry, the book Modeling Maximum Trading Profits with C++, \textit{John Wiley and Sons, Inc., Hoboken, New Jersey}, 2007, and papers in \textit{Futures Magazine} and \textit{ArXiv}.

\noindent\textit{v7f5a7@comcast.net}

\end{document}